\documentclass{amsart}

\usepackage{amsmath,amssymb,xy}

\xyoption{all}
\CompileMatrices

\theoremstyle{plain}

\newtheorem{theorem}{Theorem}[section]
\newtheorem{proposition}[theorem]{Proposition}
\newtheorem{lemma}[theorem]{Lemma}
\newtheorem{corollary}[theorem]{Corollary}
\newtheorem{conjecture}[theorem]{Conjecture}

\theoremstyle{definition}

\newtheorem{definition}[theorem]{Definition}

\newtheorem{remark}[theorem]{Remark}

\newcommand{\ST}{\SelectTips{cm}{}}

\newcommand{\A}{{\mathbf A}}
\newcommand{\F}{{\mathbf F}}
\newcommand{\kbar}{\bar{k}}
\renewcommand{\O}{{\mathcal O}}
\newcommand{\p}{{\mathfrak p}}
\newcommand{\Q}{{\mathbf Q}}
\renewcommand{\S}{{\mathcal S}}
\newcommand{\Z}{{\mathbf Z}}
\renewcommand{\o}{\underline{o}}

\newcommand{\Fp}{\F_{p}}
\newcommand{\Fpbar}{\bar{\F}_{p}}

\newcommand{\OK}{\O_{K}}
\newcommand{\Ql}{\Q_{\ell}}
\newcommand{\Qp}{\Q_{p}}
\newcommand{\Qpbar}{\bar{\Q}_{p}}
\newcommand{\Zp}{\Z_{p}}

\newcommand{\Et}{\tilde{E}}

\newcommand{\picf}{\pi^{\text{cf}}}
\newcommand{\rhob}{\bar{\rho}}
\newcommand{\univ}{\text{univ}}

\renewcommand{\th}{\text{th}}

\DeclareMathOperator{\ad}{ad}
\DeclareMathOperator{\Aut}{Aut}

\DeclareMathOperator{\Ext}{Ext}
\DeclareMathOperator{\fl}{fl}
\DeclareMathOperator{\Gal}{Gal}
\DeclareMathOperator{\GL}{GL}
\DeclareMathOperator{\rank}{rank}
\DeclareMathOperator{\Spec}{Spec}
\DeclareMathOperator{\Sym}{Sym}
\DeclareMathOperator{\tors}{tors}
\DeclareMathOperator{\ur}{ur}

\newcommand{\prank}{\rank_{p}}

\begin{document}

\title[Local torsion on elliptic curves]{Local torsion on elliptic curves and the deformation theory
of Galois representations}
\author{Chantal David and Tom Weston}
\address[Chantal David]{Department of Mathematics and Statistics, Concordia University, Montreal, QC}
\address[Tom Weston]{Department of Mathematics and Statistics, University of Massachusetts Amherst, MA}

\email[Chantal David]{cdavid@mathstat.concordia.ca}
\email[Tom Weston]{weston@math.umass.edu}

\thanks{The first author was partially supported by a NSERC Individual
Research Grant; the second author was supported by NSF grant DMS-0440708}

\maketitle

\section{Introduction}

Let $E$ be an elliptic curve over $\Q$.  Our primary goal in this paper
is to investigate for how many primes $p$ the elliptic curve $E$ possesses
a $p$-adic point of order $p$.  Simple heuristics suggest the following
conjecture.

\begin{conjecture} \label{conj:main}
Assume that $E$ does not have complex multiplication.  Fix $d \geq 1$.
Then there are finitely many primes $p$ such that there exists an
extension $K/\Qp$ of degree at most $d$ with $E(K)[p] \neq 0$.
\end{conjecture}

We expect this conjecture to be false (for sufficiently large $d$) when
$E$ does have complex multiplication.  Nevertheless, our main result is
that Conjecture~\ref{conj:main} at least holds on average.  For $A,B > 0$,
let $\S_{A,B}$ denote the set of elliptic curves with Weierstrass equations
$y^2 = x^3+ax+b$ with $a,b \in \Z$,
$|a| \leq A$ and $|b| \leq B$.  For an elliptic
curve $E$ and $x > 0$, let $\pi_{E}^{d}(x)$ denote the number of primes
$p \leq x$ such that $E$ possesses a point of order $p$ over an
extension of $\Qp$ of degree at most $d$.

\begin{theorem} \label{thm:main}
Fix $d \geq 1$.
Assume $A,B \geq x^{\frac{7}{4}+\varepsilon}$ for some $\varepsilon > 0$.
Then
$$\frac{1}{\# \S_{A,B}} \sum_{E \in \S_{A,B}} \pi^{d}_{E}(x) \ll d^2$$
as $x \to \infty$.
\end{theorem}

Our proof is essentially a precise version of the heuristic alluded
to above. We explain this heuristic for $d=1$; in this case we say
that a prime $p$ is a {\it local torsion prime} for $E$ if
$E(\Qp)[p] \neq 0$. One sees easily that for a local torsion prime
$p \geq 7$ of good reduction for $E$ one must have $a_{p}(E) = 1$.
By the Riemann hypothesis, $a_{p}(E)$ lies in an interval of length
approximately $4\sqrt{p}$, so that naively one expects $a_{p}(E)=1$
to occur with a probability of $\frac{1}{4\sqrt{p}}$.  A comparison
of the reduction exact sequence with an analysis of extension
classes suggests that when $a_{p}(E)=1$, $p$ is a local torsion prime
with probability $\frac{1}{p}$.  Combining these estimates, one
expects
$$\sum_{p} \frac{1}{4\sqrt{p}} \cdot \frac{1}{p} < \infty$$
local torsion primes for $E$, hence the conjecture.

We especially like this heuristic as each of the two parts sometimes fail.
Specifically, for some elliptic curves one has $a_{p}(E)=1$ for only
finitely many primes; in these cases, the conjecture (for $d=1$) follows
immediately.
Furthermore, for CM elliptic curves one finds that $a_{p}(E)=1$ in fact
forces $p$ to be a local torsion prime; this explains the CM exception
in the conjecture.

Although we believe that Conjecture~\ref{conj:main} is of fundamental and
independent interest in the arithmetic of elliptic curves, our
motivation for its study originated in
questions in the deformation theory of
Galois representations.
For any prime $p$ let
$$\rhob_{E,p} : \Gal(\Q_{S \cup \{p\}}/\Q) \to \GL_{2}(\Fp)$$
be the Galois representation on the $p$-torsion points of $E$; here
$\Q_{S \cup \{p\}}$ is the maximal extension of $\Q$ unramified away
from $p$ and the set $S$ of places of bad reduction for $E$.  When
$\rhob_{E,p}$ is absolutely irreducible, one can associate to
$\rhob_{E,p}$ its {\it universal deformation ring} $R_{E,p}^{\univ}$,
parameterizing all lifts of $\rhob_{E,p}$ to artinian local rings with
residue field $\Fp$.
Mazur \cite{Mazur} asked if the deformation
theory of $\rhob_{E,p}$ is {\it unobstructed} (so
that $R_{E,p}^{\univ}$ is non-canonically isomorphic to a power series
ring in three variables over $\Zp$) for all but finitely many
primes $p$.  Using work of Flach \cite{Flach}, he showed that this is
the case so long as one excludes those primes $p$ such that
$E$ possesses a point of order $p$ over a quadratic extension of $\Qp$.
Theorem~\ref{thm:main} with $d=2$ thus guarantees that this last
condition holds, on average, for only finitely many primes $p$.

Unfortunately, we can not state a deformation theoretic analogue
of Theorem~\ref{thm:main}, as the other reasons for obstructions do not
satisfy as strong of a bound.  Nevertheless, we do show that
Mazur's question has an affirmative answer for elliptic curves admitting
rational $2$-isogenies (which is basically trivial)
or $15$-isogenies (which is an amusing computation).

For another application, recall that the newform $f_{E} := \sum a_{n}(E)q^{n}$
associated to $E$ is said to possess a {\it companion form}
modulo $p$ if there is a mod $p$ eigenform $g = \sum b_{n}q^{n}$
of weight $p-1$ satisfying
$$n^{2}b_{n} \equiv na_{n} \pmod{p}$$
for all $n \geq 1$.
For $d \geq 1$,
let $\picf_{E}(x)$ denote the number of primes $p \leq x$ such
that $a_{p}(E)$ equals $\pm 1$ and
$f_{E}$ possesses a companion form modulo $p$.
Using a deep result of Gross \cite{Gross}, Theorem~\ref{thm:main}
yields the following.

\begin{theorem} \label{thm:main3}
Assume $A,B \geq x^{\frac{7}{4}+\varepsilon}$ for some $\varepsilon > 0$.
Then
$$\frac{1}{\# \S_{A,B}} \sum_{E \in \S_{A,B}} \picf_{E}(x) \ll d^2$$
as $x \to \infty$.
\end{theorem}

It would be interesting to remove the assumption that $a_{p}(E)=\pm 1$
above, although it is not immediately clear to the authors how to adapt
these methods to that question.

We present some data and simple
results on Conjecture~\ref{conj:main} in the case $d=1$ in Section~\ref{sec:data}.
In Section~\ref{sec:local}, we show that over the unramified extension
of $\Qp$ of degree $d$, residual $p$-torsion points lift to $p$-adic
$p$-torsion points with a frequency of $\frac{1}{p^{d}}$.
In Section~\ref{sec:analytic} we use this result and analytic methods
to deduce Theorem~\ref{thm:main}.
Finally, in Section~\ref{sec:defthy}
we discuss the applications to deformation theory and companion forms.

\subsection*{Acknowledgments}

The authors would like to thank Hershy Kisilevsky, Siman Wong and 
Hui June Zhu for useful discussions during the preparation of this
paper.

\section{Local torsion primes} \label{sec:data}

Fix an elliptic curve $E$ over $\Q$.  We call a prime $p$ a
{\it local torsion prime} for $E$ if $E$ possesses a point of order
$p$ over $\Qp$.  For $x > 0$ let $\pi_{E}(x)$ denote the number of
local torsion primes $p \leq x$ for $E$.
Using the Magma Computational Algebra System we computed
$\pi_{E}(10^6)$ for the 5113
elliptic curves $E$ with conductor at most $1000$.
(Note that $\pi_{E}$ is not isogeny-invariant, so that these computations
were done over isomorphism classes rather than isogeny classes.)

\vspace{0.5cm}

\begin{tabular}{c|c|cccccc}
& & \multicolumn{6}{l}{$\#E$ such that $\pi_{E}(10^6) =$} \\
Curves & \# curves & 0 & 1 & 2 & 3 & 4 & $5+$ \\ \hline
All curves & 5113 & 568 & 3687 &828 & 15 & 1 & 14 \\
Curves with no torsion & 1364 & 484 & 733 & 117 & 15 & 1 & 14
\end{tabular}

\vspace{0.5cm}

The 14 curves $E$ with $\pi_{E}(10^6) \geq 5$ were all CM curves and
had between $22$ and $36$ local torsion primes in this range.
(In fact, they all had CM by $\Q(\sqrt{-3})$, although this presumably
is an artifact of the small conductors we considered.)
Thus only one elliptic
curve without complex multiplication and conductor at most $1000$
had as many as four local torsion primes less than one million:
the elliptic curve 774D1, for which $2,3,5,7$ are local torsion primes.

We remark that the precise ratios reported here should not be taken
too seriously, as elliptic curves of small conductor are not at all
representative of all elliptic curves.  Nevertheless, this data
certainly supports Conjecture~\ref{conj:main}.  Indeed,
large local torsion primes $p$ were
quite rare: 99.1\% of the local torsion primes occurring
were $2,3,5,7$.  The only two elliptic curves in this sample with
local torsion primes $> 1000$ were 131A1 (with local torsion primes
$59$ and $4723$) and 775A1 (with local torsion prime $26993$).

We now record some simple results on local torsion primes.
Although some of the proofs rely on results in the next section,
we find it most convenient to state them here.

\begin{proposition} Let $E$ be an elliptic curve over $\Q$.
\begin{enumerate}
\item If $E(\Q)_{\tors} \neq 0$, then $E$ has finitely many local
torsion primes.  (More precisely, if $E(\Q)_{\tors} \neq 0$, then
any prime $p \geq 7$ of good
reduction for $E$ is not a local torsion prime for $E$ unless
$E(\Q)_{\tors} \cong \Z/p$.)
\item For any finite set of primes $\mathcal P$, there exists
an elliptic curve $E$
such that each element of $\mathcal P$ is
a local torsion prime for $E$.
\item If $E$ has complex multiplication, then a prime
$p > 3$ of good reduction
is a local torsion prime for $E$ if and only if $a_{p}(E)=1$.
\end{enumerate}
\end{proposition}
\begin{proof} ~
\begin{enumerate}
\item Suppose that $p \geq 7$ is a local torsion prime of $E$ at
which $E$ has good reduction. Since $p \geq 3$, the formal group of
$E$ over $\Zp$ is torsion-free, so that the natural map
$$E(\Qp)_{\tors} \to E(\Fp)$$ is injective.  In particular, $p$
divides $\#E(\Fp)$; since $p \geq 7$, it now follows from the Riemann
hypothesis that $E(\Fp)$ must have exactly $p$ elements. Since
$E(\Q)_{\tors}$ is non-trivial and injects into $E(\Fp)$, we
conclude that $E(\Q)_{\tors} \cong \Z/p$, as desired.
\item For each $p \in {\mathcal P}$, fix $a_p,b_p \in \Fp$ so that
the elliptic curve $E_{a_p,b_p}$ with Weierstrass equation
$y^2=x^3+a_px+b_p$ has precisely $p$ points over $\Fp$.  (The
existence of such a curve is guaranteed by Deuring's Theorem which
gives an exact formula for the number of elliptic curves over $\Fp$
with $p+1-r$ points for $|r| \leq 2 \sqrt{p}$; see for example
\cite{Lenstra} for a statement of Deuring's Theorem.) Let
$(A_{p},B_{p}) \in \Z/p^2 \times \Z/p^2$ be one of the pairs lifting
$(a_{p},b_{p})$ in Corollary~\ref{cor:alg}.  Let $A,B \in \Z$ be
such that $A \equiv A_{p} \pmod{p^{2}}$, $B \equiv B_{p}
\pmod{p^{2}}$ for all $p \in {\mathcal P}$.  Then by
Corollary~\ref{cor:alg} each $p \in {\mathcal P}$ is a local torsion
prime for $E_{A,B}$.
\item This is clear from Lemma~\ref{lemma:canonical_lift} as
a CM elliptic curve is a canonical lift of each of its reductions.
Alternately, it follows from the fact that for an ordinary prime $p$
the mod $p$ Galois representation of a CM elliptic curve
is abelian when restricted to a decomposition group at $p$.
\end{enumerate}
\end{proof}

\section{Local torsion on elliptic curves} \label{sec:local}

Fix a prime $p > 3$
and a finite extension $k$ of $\Fp$ of degree $d$.
Let $W$ denote the ring of Witt vectors over $k$; we write
$K$ for the field of fractions of $W$.  Note that $k = W/pW$;
we set $W_{2} = W/p^2W$.

Let $E$ be an elliptic curve over $W$; that is, $E$ is an elliptic
curve over $K$ with good reduction.  The next lemma gives a criterion
for $E$ to possess a $K$-rational point of order $p$.  For a finite
abelian group $M$ we write $\prank M$ for the $p$-rank of $M$
(that is, for the $\Fp$-dimension of $M \otimes_{\Z} \Fp$ or equivalently
of the $p$-torsion subgroup $M[p]$).

\begin{lemma} \label{lemma:prank}
Let $E$ be an elliptic curve over $W$.  Then:
$$\prank E(W_2) = \begin{cases}
d & \mbox{if $E(K)[p] = 0$;} \\
d+1 & \mbox{if $E(K)[p] \neq 0$.} \end{cases}$$
\end{lemma}
\begin{proof}
Consider the diagram (with exact rows)
$$\ST\xymatrix{
{0} \ar[r] & {E_{0}(pW)} \ar[r] \ar[d] & {E(K)} \ar[r] \ar[d] & {E(k)} \ar[r] \ar@{=}[d] & {0} \\
{0} \ar[r] & {E_{0}(pW/p^2)} \ar[r] & {E(W_2)} \ar[r] & {E(k)} \ar[r] & {0}}$$
with $E_{0}$ the formal group of $E$ over $\OK$.  Since $W$ is unramified
over $\Zp$, it follows from \cite[Theorem IV.6.4b]{Silverman} that the
reduction map
$E_{0}(pW) \to E_{0}(pW/p^2)$ can be identified with the natural map
$\Zp^d \to (\Z/p)^d$.
In particular, taking $p$-torsion and applying the snake lemma
we obtain a commutative diagram
$$\ST\xymatrix{
{0} \ar[r] & {0} \ar[r] \ar[d] & {E(K)[p]} \ar[r] \ar[d] & {E(k)[p]} \ar[r] \ar@{=}[d] & {(\Z/p)^{d}}\ar@{=}[d] \\
{0} \ar[r] & {(\Z/p)^d} \ar[r] & {E(W_2)[p]} \ar[r] & {E(k)[p]} \ar[r] & {(\Z/p)^{d}}}$$
It follows that
$$\prank E(W_2) = d + \prank E(K)[p].$$
Since the field $K(E[p])$ contains $K(\mu_{p})$ and thus
is ramified over $\Qp$,
$E(K)[p]$ it at most one-dimensional and the lemma follows.
\end{proof}

We fix now an elliptic curve $E$ over $k$ such that
$E(k)[p] \neq 0$.  We first determine how many lifts of
$E$ to $W_{2}$ have $p$-rank $d+1$.

\begin{lemma} \label{lemma:canonical_lift}
Let $E$ be an elliptic curve over $k$ such
that $E(k)$ has a non-trivial $p$-torsion point.  Then precisely one
 of the $p^d$ isomorphism classes of
lifts of $E$ to an elliptic curve $\Et$ over $W_2$ satisfies
$\prank \Et(W_2) = d+1$.
\end{lemma}
\begin{proof}
Since $E$ is necessarily ordinary, there is a canonical exact
sequence of finite flat group schemes
$$0 \to \mu_{p} \to E[p] \to \Z/p \to 0$$
over $\kbar$.  In fact, since $E(k)$ has a point of order $p$, it
follows that this exact sequence exists already over $k$ and is split;
that is, $E[p] \cong \mu_{p} \times \Z/p$ as
$k$-group schemes.

By the Serre--Tate theorem (see \cite[Theorem 1.2.1]{Katz}, for example),
the lifts of $E$ to elliptic curves over $W_2$ are
parameterized by the lifts of the $p$-divisible group
$E[p^{\infty}]$ to a $p$-divisible group over $W_2$.  Any such lift is
an extension of $\Qp/\Zp$ by $\mu_{p^{\infty}}$, so that these lifts
correspond to elements of
$$\underset{n}{\varprojlim} \Ext^{1}(\Z/p^{n},\mu_{p^{n}})$$
where the extensions are computed in the category of finite flat
group schemes over $W_2$.  It is equivalent to compute these extensions
under the flat topology, so that
$$\Ext^{1}(\Z/p^{n},\mu_{p^{n}}) \cong H^{1}_{\fl}(\Spec W_2,
\mu_{p^{n}}) \cong W_{2}^{\times}/W_{2}^{\times p^{n}} \cong p \cdot W_2
\cong (\Z/p)^d.$$
Therefore the natural map
\begin{equation} \label{eq:ext}
\underset{n}{\varprojlim} \Ext^{1}(\Z/p^{n},\mu_{p^{n}})
\to \Ext^{1}(\Z/p,\mu_{p})
\end{equation}
is an isomorphism, and these groups have order $p^d$.
In particular, this proves that
$E$ has precisely $p^d$ isomorphism classes of lifts to $W_{2}$.

Let $\Et$ be a lift of $E$ to $W_2$.  One sees immediately from the
proof of Lemma~\ref{lemma:prank} that $\Et(W_2)$ has $p$-rank $d+1$
if and only if the canonical exact sequence
$$0 \to \mu_{p} \to \Et[p] \to \Z/p \to 0$$
splits.  By (\ref{eq:ext}), this occurs for precisely one lift $\Et$,
as claimed.  (This lift is usually called the {\it canonical lift} of
$E$ to $W_2$.)
\end{proof}

\begin{remark}
We note that it follows from the above proof that the isomorphism class
of an elliptic curve over $W_2$ is determined by its $j$-invariant and
the isomorphism class of its reduction to $k$.  This fact can certainly
be proven in for more elementary ways.
\end{remark}

For $a,b$ elements of some ring $R$, we write $E_{a,b}$ for the
projective plane curve with affine Weierstrass equation $y^{2} =
x^{3} + ax + b$.  The next lemma applies Lemma~\ref{lemma:canonical_lift}
to obtain a lifting result for these cubics.

\begin{proposition} \label{prop:local}
Let $E_{a,b}$ be an elliptic curve over $k$ such that
$j(E_{a,b}) \neq 0, 1728$.  Assume that $E_{a,b}(k)$
has an element of order $p$.
Then there are exactly $p^{d}$
distinct pairs $(A_{i},B_{i}) \in W_2 \times W_2$
such that $(A_{i},B_{i}) \equiv (a,b) \mod{p}$ and
$$\rank_{p} E_{A_{i},B_{i}}(W_2) = d+1.$$
\end{proposition}
\begin{proof}
Consider the closed subscheme $S$ of $\A^3_{W_2}$ (with coordinates
$j,A,B$) defined by the
$j$-invariant equation
$$j = 6912 \cdot \frac{A^3}{4A^3+27B^2}.$$
The surface $S \times_{W_2} k$ is smooth away from $j=0,1728$,
so that it follows from Hensel's lemma that for each of the $p^d$
lifts $\tilde{j}$ of $j(E_{a,b})$ to $W_2$, there are precisely $p^d$ lifts of
$(a,b)$ to $W_2 \times W_2$
such that the corresponding elliptic curve over $W_2$ has
$j$-invariant $\tilde{j}$.  In particular, each of the $p^d$
possible $j$-invariants lifting $j(E)$ to $W_2$ occurs for some $E_{A,B}$
with $(A,B) \in W_{2} \times W_{2}$.
Since by Lemma~\ref{lemma:canonical_lift} there
are precisely $p^d$ isomorphism classes of lifts of $E_{a,b}$ to $W_{2}$,
it follows
that the isomorphism class of a lift of $E$ is determined by its
$j$-invariant.  (Alternately, since the twists  of an elliptic curve $E$
over a ring $R$ are parameterized by $H^{1}_{\fl}(\Spec R,\Aut E)$ and
one has $p \nmid \# \Aut E_{a,b}$ since $p \geq 5$, one
can deduce this from the fact that
the categories
of \'etale sheaves on $\Spec k$ and $\Spec W$ are equivalent.)
In particular, each isomorphism class occurs for precisely $p^d$ values
of $A,B$.
By Lemma~\ref{lemma:canonical_lift} precisely one
of these isomorphism classes has $p$-rank $d+1$, which yields the
proposition.
\end{proof}

For our analytic arguments we will be interested only in elliptic
curves over $\Zp$.  In this case Proposition~\ref{prop:local} yields
the following.

\begin{corollary} \label{cor:alg}
Let $E_{a,b}$ be an ordinary
elliptic curve over $\Fp$ such that $j(E_{a,b}) \neq 0,1728$.
Let $k$ be an extension of $\Fp$ of degree $d$ such that
$E(k)[p] \neq 0$; set $W_2 = W/p^2$ with $W$ the ring of Witt vectors
of $k$.  Then
there are exactly $p$ distinct pairs $(A_{i},B_{i}) \in \Z/p^2 \times \Z/p^2$
such that $(A_{i},B_{i}) \equiv (a,b) \mod{p}$ and
$$\rank_{p} E_{A_{i},B_{i}}(W_2) = d+1.$$
\end{corollary}
\begin{proof}
The canonical lift of $E_{a,b}$ to $W_{2}$ has $j$-invariant in
$\Z/p^2$.
Thus by Hensel's lemma as before
there are precisely $p$ pairs lying in $\Z/p^2 \times \Z/p^2$
among the $p^d$ pairs of
Proposition~\ref{prop:local}.  The corollary follows.
\end{proof}

\section{Analytic arguments} \label{sec:analytic}

Fix $d \geq 1$ and as before
let $W_{2}$ denote the Witt vectors of length two over
a finite field $k$ of order $p^d$.

\begin{definition}
We write $\nu_{p}(d)$ for the number of pairs
$(a,b) \in \Z/p^2 \times \Z/p^2$ such that $E_{a,b}$ is an elliptic
curve with $\rank_{p} E_{a,b}(W_2) = d+1$.
\end{definition}

The following estimates will be used
in our main result.

\begin{lemma} \label{lemma:estimates}
We have
$$\sum_{p \leq x} \nu_{d}(p) \ll dx^{7/2}\log x;$$
$$\sum_{p \leq x} \frac{\nu_{d}(p)}{p^2} \ll dx^{3/2}\log x;$$
$$\sum_{p \leq x} \frac{\nu_{d}(p)}{p^4} \ll d.$$
\end{lemma}
\begin{proof}
Write $\nu_{d}(p) = \nu'_{d}(p) + \nu^{0}_{d}(p) + \nu^{1728}_{d}(p)$
where $\nu'_{d}(p)$ (resp.\ $\nu^{0}_{d}(p)$, resp.\ $\nu^{1728}_{d}(p)$)
denotes the number of pairs $(a,b) \in \Z/p^2 \times \Z/p^2$ such that
$E_{a,b}(W_2)$ has $p$-rank $d+1$ and $E_{a,b}$ does not have
$j$-invariant $0$ or $1728$ (resp.\ has $j$-invariant $0$,
resp.\ has $j$-invariant $1728$).
It suffices to prove the lemma separately for each of these three
functions.

We begin with $\nu'_{d}(p)$.  Since an elliptic curve $E$ over $\Fp$
has a point of order $p$ over $\F_{p^{d}}$ if and only if $a_p(E)^d
\equiv 1 \mod p$ (see Lemma~\ref{lemma:apd} below), by
Corollary~\ref{cor:alg} we have
\begin{align*}
\nu'_{d}(p) &= p \cdot \# \bigl\{ (a,b) \in \Fp^{\times} \times \Fp^{\times} \,;\,
a_p(E_{a,b})^d \equiv 1 \bmod{p} \bigr\} \\
&= p \cdot \sum_{\substack{|r| < 2\sqrt{p} \\ r^d \equiv 1 \bmod{p}}} (p-1)H(r^2-4p)\\
\end{align*}
by Deuring's theorem; see for example \cite{Lenstra}.
(Note that restricting to $(a,b) \in \Fp^{\times} \times \Fp^{\times}$
has eliminated those elliptic curves with $j$-invariant $0$ or $1728$.)
The class number formula together with the trivial bound
$L(1,\chi_{D}) \ll \log D$ yields
$H(r^2-4p) \ll \sqrt{p}\log^{2}p$.  As (for large enough $p$) there are
at most $d$ such $r$, it follows that
$$\nu'_{d}(p) \ll dp^{5/2}\log^{2}p.$$
The asserted bounds for $\nu_{d}'(p)$ all follow easily from this.

Consider now $\nu^{0}_{d}(p)$.  An elliptic curve $E$ over $\Z/p^2$
with $j$-invariant $0$ is always the canonical lift
of its reduction, so that by the proof of Lemma~\ref{lemma:canonical_lift}
it has $p$-rank $d+1$ over $W$
precisely when
$$a_{p}(E \times_{\Z/p^2} \Fp)^{d} \equiv 1 \pmod{p}.$$
(Note that in this case Corollary~\ref{cor:alg} does not apply.)
Since a curve $E_{a,b}$ over $\Fp$ has $j$-invariant $0$ if and only if
$b =0$ and
each such curve over $\Fp$ has $p^2$ lifts to $\Z/p^2$, we obtain
$$\nu^{0}_{d}(p) = p^2 \cdot \sum_{
\substack{|r| < 2\sqrt{p} \\ r^d \equiv 1 \bmod{p}}}
\# \{ a \in \Fp^{\times} ; a_p(E_{a,0}) = r\}.$$
Thus
\begin{align*}
\sum_{p \leq x} \frac{\nu^{0}_{d}(p)}{p^2} &=
\sum_{p \leq x} \sum_{\substack{|r| < 2\sqrt{p}\\r^d \equiv 1 \bmod{p}}} \sum_{a \in \Fp^{\times}}
\chi_{E_{a,0}}^{r}(p)
\end{align*}
where $\chi_{E}^{r}(p)$ is $1$ if
$a_{p}(E) = r$ and $0$ otherwise.
Reversing the order of summation, we obtain
\begin{align*}
\sum_{p \leq x} \frac{\nu^{0}_{d}(p)}{p^2} &=
\sum_{a \leq x} \sum_{a \leq p \leq x} \sum_{
\substack{|r| < 2\sqrt{p}\\ r^d \equiv 1 \bmod{p}}}
\chi_{E_{a,0}}^{r}(p) \\
&\ll d \cdot \underset{r}{\max} \left( \sum_{a \leq x} \# \bigl\{
p \leq x \,;\,
a_{p}(E_{a,0}) = r \bigr\} \right).
\end{align*}
An elementary argument using the description of $a_p(E_{a,0})$
in terms of the quadratic form $X^2+Y^2$ shows that
$$\# \bigl\{ p \leq x ; a_{p}(E_{a,0}) = r \bigr\} \ll \sqrt{x}$$
so that we obtain
\begin{align*}
\sum_{p \leq x} \frac{\nu^{0}_{d}(p)}{p^2} &\ll
d \cdot \sum_{a \leq x} \sqrt{x} \\
&\ll dx^{3/2}.
\end{align*}
Applying partial summation one now easily obtains
$$\sum_{\p \leq x} \frac{\nu^{0}_{d}(p)}{p^4} \ll d$$
$$\sum_{\p \leq x} \nu^{0}_{d}(p) \ll dx^{7/2}$$
which suffice for the lemma.
An entirely similar argument gives the same bounds
for $\nu_{d}^{1728}$.
\end{proof}

The next lemma was used in the preceding proof.  There are of course
many proofs; we give one in the spirit of Section 3. 
Alternately, one could
use the Weil conjectures and look at the $p$-divisibility of $\#
E(\F_{p^d}) = p^d+1-(\alpha_p^d+\beta_p^d)$, where $\alpha_p,
\beta_p$ are the roots of the characteristic polynomial
$x^2-a_p(E)x+p$ over $\Fp$.

\begin{lemma} \label{lemma:apd}
Let $E$ be an elliptic curve over $\Fp$.  Then
$E(\F_{p^d})[p] \neq 0$ if and only if $a_p(E)^d \equiv 1 \mod p$.
\end{lemma}
\begin{proof}
If $a_p(E)=0$, then $E$ is supersingular and never has $p$-torsion over
any extension of $\Fp$, so we may assume that $a_p(E) \neq 0$.  In this case,
the $p$-torsion of $E$ decomposes as
$$E[p] = \Z/p(\chi) \oplus \mu_p(\chi^{-1})$$
where $\chi$ is the character of $\Gal(\Fpbar/\Fp)$ sending a Frobenius
to $a_p(E)$.  (See \cite{Gross}, for example.)
Now $E(\F_{p^d})[p] \neq 0$ if and only if $\chi$ is trivial on
$\Gal(\Fpbar/\F_{p^d})$; since this Galois group is generated by
the $d^{\th}$ power of Frobenius, the lemma follows.
\end{proof}

Recall that for an elliptic curve $E$,
$\pi_{E}^{d}(x)$ denotes the number of primes $p \leq x$ such
that $E$ possesses a point of order $p$ over an extension
of $\Qp$ of degree at most $d$.  Let $\pi_{E}^{d,\ur}(x)$
denote the number of primes $p \leq x$ such that $E$ possesses a
point of order $p$ over an unramified extension of $\Qp$ of degree
at most $d$.

\begin{lemma} \label{lemma:ram}
Let $E$ be an elliptic curve.  Then
$$\pi_{E}^{d}(x) - \pi_{E}^{d,\ur}(x)$$
is bounded by the number of primes $p \leq d$.
\end{lemma}
\begin{proof}
Fix a prime $p > d$ and a ramified extension $K/\Qp$ of degree
at most $d$.  To prove the lemma it suffices to show
that if $E(K)[p] \neq 0$, then already $E(K^{\ur})[p] \neq 0$
with $K^{\ur}$ the maximal unramified subfield of $K$.
For this, recall that the restriction of the
$p$-torsion Galois representation of $E$ to a decomposition group at $p$
has the form
$$\left( \begin{array}{cc}
\varepsilon\chi & * \\ 0 & \chi^{-1} \end{array}
\right)$$
with $\varepsilon$ the cyclotomic character, $\chi$  an unramified
character and $*$ either trivial or wildly ramified.  $K$ is not wildly
ramified (since $d<p$) and does not contain $\mu_{p}$, so that
the only way that $E(K)[p]$ can be non-zero is if $*$ is trivial
and $\chi$ factors through $\Gal(K/\Qp)$.  Since
$\chi$ is unramified, it must thus factor through $\Gal(K^{\ur}/\Qp)$
as well, as desired.
\end{proof}

Recall that $\S_{A,B}$ denotes the set of elliptic curves
$E_{a,b}$ with $a,b \in \Z$ and $|a| \leq A$, $|b| \leq B$.

\begin{proposition}
Fix $A,B > 0$.  Then
$$\frac{1}{\# \S_{A,B}} \sum_{E \in \S_{A,B}} \pi_{E}^{d}(x) =
O(d^2) + d \cdot
O\left( \left(\frac{1}{A} + \frac{1}{B} \right)x^{3/2} \log x +
\frac{1}{AB} x^{7/2} \log x \right).$$
\end{proposition}
\begin{proof}
By Lemma~\ref{lemma:ram} we have
$$\frac{1}{\# \S_{A,B}} \sum_{E \in \S_{A,B}} \pi_{E}^{d}(x) =
\frac{1}{\# \S_{A,B}} \sum_{E \in \S_{A,B}} \pi_{E}^{d,\ur}(x) + d \cdot
O(1).$$
It thus suffices to consider the latter sum.

Let $K_p^d$ denote the unramified extension of $\Qp$ of degree $d$ and
let $\tilde{\pi}_{E}^{d}(x)$ denote the number of $p \leq x$ such
that $E(K_p^d)[p] \neq 0$.
By Lemma 3.1, $E(K_p^d)[p] \neq 0$ can be detected by looking at $E$
over $\Z/p^2$.  Since there are
$$\left( \frac{2A}{p^{2}} + O(1) \right) \cdot
\left( \frac{2B}{p^{2}} + O(1) \right)$$
elliptic curves in $\S_{A,B}$ reducing to any fixed elliptic curve
over $\Z/p^2$, applying Lemma~\ref{lemma:prank} and reversing the order of
summation we find that
\begin{align*}
\frac{1}{\# \S_{A,B}} \sum_{E \in \S_{A,B}} \tilde{\pi}_{E}^{d}(x) &=
\frac{1}{\# \S_{A,B}} \sum_{p \leq x} \# \bigl\{ E \in \S_{A,B} \,;\,
E(K_p^d)[p] \neq 0 \bigr\} \\
&= \frac{1}{\# \S_{A,B}} \sum_{p \leq x}
\left( \frac{2A}{p^{2}} + O(1) \right) \cdot
\left( \frac{2B}{p^{2}} + O(1) \right) \nu_{d}(p) \\
&= \sum_{p \leq x} \frac{\nu_{d}(p)}{p^{4}} +
O \left( \left(\frac{1}{A}+\frac{1}{B}\right)
\sum_{p \leq x} \frac{\nu_{d}(p)}{p^{2}} + \frac{1}{AB}
\sum_{p \leq x} \nu_{d}(p) \right)
\end{align*}
since
$$\# \S_{A,B} = 4AB(1+ \o(1)).$$
Applying the estimates of Lemma~\ref{lemma:estimates}
and summing over $d$ now yields the result.
\end{proof}

\begin{corollary}
If $A,B \gg x^{7/4+\varepsilon}$, then
$$
\frac{1}{\# \S_{A,B}} \sum_{E \in \S_{A,B}} \pi_{E}^{d}(x) \ll
d^2$$
as $x \to \infty$.
\end{corollary}

\section{Deformation theory} \label{sec:defthy}

Fix an elliptic curve $E$ over $\Q$ without complex multiplication.
Let $S$ denote the set of primes at which $E$ has bad reduction.
For a prime $p$ consider the Galois representation
$$\rhob_{E,p} : \Gal(\Q_{S \cup \{p\}}/\Q) \to \GL_{2}(\Fp)$$
giving the Galois action on the $p$-torsion points $E[p]$.
A {\it deformation} of $\rhob_{E,p}$
to a local ring $R$ with residue field $\Fp$ is a
strict equivalence class of Galois representations
$$\Gal(\Q_{S \cup \{p\}}/\Q) \to \GL_{2}(R)$$
which yield $\rhob_{E,p}$ on composition with the reduction map
$$\GL_{2}(R) \to \GL_{2}(\Fp);$$
here two Galois representations are considered to be
strictly equivalent if they are
conjugate by a matrix which reduces to the identity in $\GL_2(\Fp)$.

When $\rhob_{E,p}$ is absolutely irreducible (which is true for
sufficiently large $p$) it is known \cite[Section 1.2]{Mazur2}
that there exists a {\it universal
deformation}
$$\rho_{E,p}^{\univ} : \Gal(\Q_{S \cup \{p\}}/\Q) \to \GL_{2}(R_{E,p}^{\univ})$$
in the category of inverse limits of artinian
local rings with residue field $\Fp$; the ring $R_{E,p}^{\univ}$ is called the
{\it universal deformation ring} of $\rhob_{E,p}$.
Thus any deformation of
$\rhob_{E,p}$ to such a ring $R$ is obtained
from $\rho^{\univ}$ by composition with a unique map $R_{E,p}^{\univ} \to R$.

It is of fundamental interest to understand the structure of
$R_{E,p}^{\univ}$.  Let $\ad \rhob_{E,p}$ denote the adjoint
representation of $\rhob_{E,p}$ and set
$$d_{i} := \dim_{\Fp} H^{i}(\Q_{S \cup \{p\}}/\Q,\ad \rhob_{E,p}).$$
We have the following result of Mazur \cite[Section 1.6 and Section 1.10]{Mazur2}.

\begin{proposition}
With notation as above,
$R_{E,p}^{\univ}$ is a quotient of a power series ring in
$d_{1}$ variables over $\Zp$ by an ideal generated by at most
$d_{2}$ elements.  Furthermore, one has $d_{1}-d_{2}=3$.  In particular, if
$d_{2}=0$ (in which case one says that the deformation theory of
$\rhob_{E,p}$ is unobstructed), then $R_{E,p}^{\univ}$ is
(non-canonically) isomorphic to a power series ring in three variables
over $\Zp$.
\end{proposition}

In \cite{Mazur}, Mazur further asks if the deformation theory of
$\rhob_{E,p}$ is unobstructed
for all but finitely many primes $p$.
(We remark that the analogous question for modular forms of weight
at least $3$ was answered affirmatively in \cite{TWdef}.)
Using Poitou--Tate duality and work of Flach \cite{Flach}
on symmetric square Selmer
groups, he proved the following
result.

\begin{proposition} \label{prop:mazur}
Let $p$ be a prime of good reduction for $E$ such that
$\rhob_{E,p}$ is absolutely irreducible.
If the deformation theory of $\rhob_{E,p}$ is obstructed, then one of
the following holds:
\begin{enumerate}
\item $p \leq 3$;
\item $p \in S$;
\item $\ell \equiv 1 \pmod{p}$ for some $\ell \in S$;
\item The Galois representation $\rhob_{E,p}$ is not surjective;
\item $p$ divides the (appropriately normalized) special value
$L(\Sym^2 E,2)/\Omega$ of the symmetric square $L$-function of $E$;
\item $H^{0}(\Ql,\Sym^2 E[p]) \neq 0$ for some $\ell \in S \cup \{p\}$.
\end{enumerate}
Furthermore each of these conditions holds for all but finitely many
primes $p$, with the possible exception of the vanishing of
$H^{0}(\Qp,\Sym^{2} E[p])$ in (6).
\end{proposition}

To answer Mazur's question it thus suffices to show that
$H^{0}(\Qp,\Sym^2 E[p])=0$ for all but finitely many $p$.
The next proposition gives some reformulations of this condition.

\begin{proposition} \label{prop:mazurgross}
Let $p$ be an odd prime of good reduction for $E$.  Then
the following are equivalent:
\begin{enumerate}
\item $H^{0}(\Qp,\Sym^2 E[p]) \neq 0$;
\item $a_{p}(E) = \pm 1$ and the restriction of $\rhob_{E,p}$ to
a decomposition group at $p$ is semi-simple;
\item $a_{p}(E) = \pm 1$ and $f_{E}$ possesses a companion form modulo
$p$.
\end{enumerate}
If $p \geq 7$, then these conditions are further equivalent to:
\begin{enumerate}
\setcounter{enumi}{3}
\item $E(K)[p] \neq 0$ with $K$ some quadratic extension of
$\Qp$.
\end{enumerate}
\end{proposition}
\begin{proof}
The equivalence of (1) and (2) is \cite[Lemma of p.\ 172]{Mazur}.
The equivalence of (2) and (3) is the main result of
\cite{Gross}.  Finally, for the equivalence of (1) and (4) we modify
the argument  of \cite[Lemma 2.3(i)]{coates}.
It is clear that if $x \in E(K)[p]$ is non-trivial, then
$x \otimes x$ gives a non-zero element of
$H^{0}(\Qp,\Sym^2 E[p])$.  For the converse, suppose that
$t \in H^{0}(\Qp,\Sym^2 E[p])$ is non-trivial.  Let $\varepsilon$ denote the
mod $p$ cyclotomic character.  Since $p \geq 7$ we may
choose $\sigma \in \Gal(\Qpbar/\Qp)$
such that $\varepsilon(\sigma)^4 \neq 1$.  Let
$\lambda,\mu \in \Fpbar$ be the eigenvalues of $\sigma$ acting
on $E[p]$.  Then the eigenvalues of $\sigma$ acting on
$\Sym^{2}E[p]$ are $\lambda^2$, $\lambda\mu = \varepsilon(\sigma)$,
and $\mu^2$.
Since $t$ is fixed by $\sigma$, one of these must equal $1$; as
$\varepsilon(\sigma) \neq 1$ we may assume without loss of generality
that $\lambda^2 = 1$.

Now consider $\sigma^2$, which has $E[p]$ eigenvalues
$\lambda^2=1$ and $\varepsilon(\sigma)^2$.  These eigenvalues
are distinct, so we may choose a basis $x,y$ of $E[p]$ with
$\sigma^2$-eigenvalues $1$ and $\varepsilon(\sigma)^2$, respectively.
Since $\varepsilon(\sigma)^4 \neq 1$, one computes easily that
the $\sigma^2$-invariant subspace of $\Sym^2 E[p]$ is spanned by
$x \otimes x$.  Since $t$ must lie in this subspace, it follows that
$x \otimes x$ is fixed by all of $\Gal(\Qpbar/\Qp)$.  Thus
$\Gal(\Qpbar/\Qp)$ acts on $x$ by a quadratic character, so that
$x$ is defined over a quadratic extension of $\Qp$.
\end{proof}

Combined with Proposition~\ref{prop:mazurgross},
Theorem~\ref{thm:main} immediately yields Theorem~\ref{thm:main3}.
Furthermore, it shows that the exceptional case
in Proposition~\ref{prop:mazur} occurs on average finitely often.
Unfortunately, already the condition (2) of
Proposition~\ref{prop:mazur} fails to satisfy such a strong bound,
so that we can not state such a result for the deformation theory of $E$.

Nevertheless, for certain elliptic curves one has $a_{p}(E)=\pm 1$
for finitely many $p$, so that
one can verify Conjecture~\ref{conj:main} for
$d=2$ for these curves and thus settle these questions.
More precisely, we have the following.
(For an alternate proof of the $15$-isogeny case, see
\cite[Theorem 1.4]{David}.)

\begin{proposition} \label{prop:torsion}
Let $E$ be an elliptic curve over $\Q$ without complex multiplication.
Assume that $E$ admits either a rational $2$-isogeny or a rational
$15$-isogeny.  Then
Conjecture~\ref{conj:main} holds for $E$ with $d=2$.  In particular,
the deformation theory of $\rhob_{E,p}$ is unobstructed for all but
finitely many primes $p$.
\end{proposition}
\begin{proof}
We in fact show that $E$ does not have a point of order $p$
over $\F_{p^2}$ for all but finitely many $p$.
Since for $p \geq 7$
this occurs if and only if $a_{p}(E) = \pm 1$, it suffices
to show that occurs for only finitely many primes $p$.

If $E$ possesses a rational $2$-isogeny, then it has rational
$2$-torsion.  Thus in particular
$E(\Fp)[2] \neq 0$ for all but finitely many primes $p$.
Fix such a $p > 2$.  Then
$$a_{p}(E) = p + 1 - \#E(\Fp) \equiv 0 \pmod{2},$$
so that certainly $a_{p}(E) \neq \pm 1$, as desired.

Now let $E$ be an elliptic curve admitting a $15$-isogeny.
The modular curve $X_0(15)$ has four non-cuspidal rational points; as the
four elliptic curves in the isogeny class 50A of \cite{Cremona} admit rational
$15$-isogenies and are distinct over $\bar{\Q}$, it follows that every
elliptic curve admitting a rational $15$-isogeny is isogenous to a twist
of the curve 50A1.  Since isogeny does not change Fourier coefficients
and twisting changes them only up to sign (and introducing some zeroes),
it suffices to prove the proposition for the elliptic
curve 50A1.

Let $E$ denote the elliptic curve $50A1$; it
has Weierstrass equation
$$y^{2} + xy + y = x^3 - x - 2.$$
Then $E$ has the rational $3$-torsion points
$$\{ (2,1), (2,-4) \} \in E[3]$$
and admits a rational $5$-isogeny.
Thus
$$\rhob_{E,3} \cong \left( \begin{array}{cc} 1 & * \\ 0 & \varepsilon
\end{array} \right) \qquad
\rhob_{E,5} \cong \left( \begin{array}{cc} \chi & * \\ 0 & \varepsilon\chi^{-1}
\end{array} \right)$$
with $\varepsilon$ cyclotomic and $\chi$ some Dirichlet character.
One finds that $E$ has a $5$-torsion point
$$P = (-\zeta_5^3 - \zeta_5^2 - 1, -\zeta_5^3 - 2\zeta_5-1)$$
with $\zeta_5$ a primitive $5^\text{th}$ root of unity.
Using this one explicitly computes $\chi$ as a Dirichlet character
of conductor $5$; combined with
the knowledge of $\rhob_{E,3}$, we conclude that
$$\rhob_{E,15} \cong \left( \begin{array}{cc} \chi_{1} & * \\ 0 &
\chi_{2}  \end{array} \right)$$
with $\chi_{1},\chi_{2}$ Dirichlet characters of conductor $15$ given by
\begin{center}
\begin{tabular}{c|cccccccc}
$n$ & 1 & 2 & 4 & 7 & 8 & 11 & 13 & 14 \\ \hline
$\chi_{1}(n)$ & 1 & 13 & 4 & 13 & 7 & 1 & 7 & 4 \\
$\chi_{2}(n)$ & 1 & 14 & 1 & 4 & 14 & 11 & 4 & 11
\end{tabular}
\end{center}
In particular, we have
$$a_{p}(E) \equiv \chi_{1}(p) + \chi_{2}(p) \pmod{15}$$
for all $p \geq 7$.  Thus
$$a_{p}(E) \in \{0,2,5,6,11,12 \}$$
modulo $15$, so that $a_{p}(E) \neq \pm 1$ for all $p \geq 7$, as desired.
\end{proof}

\end{document}